\documentclass[12pt]{amsart}
\usepackage{amsmath,amscd,amssymb,amsfonts}
\usepackage[hypertex, pagebackref=true]{hyperref}
\hypersetup{colorlinks=false}

\newcommand{\bA}{{\mathbb A}}
\newcommand{\bC}{{\mathbb C}}

\def\bF{\mathbb F}

\newcommand{\bN}{{\mathbb N}}
\newcommand{\bP}{{\mathbb P}}
\newcommand{\bR}{{\mathbb R}}
\newcommand{\bQ}{{\mathbb Q}}
\newcommand{\bZ}{{\mathbb Z}}

\newcommand{\cD}{{\mathcal D}}

\newcommand{\cJ}{{\mathcal J}}
\newcommand{\cO}{{\mathcal O}}

\newcommand{\cS}{{\mathcal S}}

\newcommand{\Hom}{{\rm{Hom}}}

\newcommand{\wti}{\widetilde}

\newcommand{\Gr}{\text{\rm Gr}}

\newcommand{\Sp}{{Sp}}

\newcommand{\codim}{\hbox{\rm codim}\,}

\newcommand{\ra}{\rightarrow}
\newcommand{\rndown}[1] {\lfloor {#1} \rfloor}
\newcommand{\rndup}[1] {\lceil {#1} \rceil}

\def\lct{{{lct}}\,}
\def\fpt{{{fpt}}\,}
\def\pa{\partial}
\def\al{\alpha}
\def\cX{\mathcal X}
\def\ti{\tilde}

\def\eps{\epsilon}
\def\lam{\lambda}
\def\me{\medskip}
\def\ni{\noindent}

\title[Singularity invariants related to Milnor fibers: survey]{Singularity invariants related to Milnor fibers: survey}
\author{Nero Budur}
\address{Department of Mathematics,
University of Notre Dame, 255 Hurley Hall, IN 46556, USA} \email{nbudur@nd.edu}

\thanks {This work was partially supported by NSF and NSA}

\begin{document}

\begin{abstract}
 This brief survey of some singularity invariants related to Milnor fibers should serve as a quick guide to references. We attempt to place things into a wide geometric context while leaving technicalities aside. We focus on relations among different invariants and on the practical aspect of computing them.
\end{abstract}

\maketitle


\tableofcontents


	Trivia: in how many different ways can the log canonical threshold of a polynomial be computed ? At least  6 ways in general, plus 4 more ways with some luck.	
	
	\medskip
	Singularity theory is a subject deeply connected with many other fields of mathematics. We give a brief 		survey of some singularity invariants related to Milnor fibers that should serve as a quick guide to references. This is by no means an exhaustive survey and many topics are left out. What we offer in this survey is an attempt to place things into a wide geometric context while leaving technicalities aside. We focus on relations among different invariants and on the practical aspect of computing them.	Along the way we recall some questions  to serve as food for thought.

	\me
	To achieve the goals we set with this survey, we pay a price. This is not a historical survey, in the sense that general references are mentioned when available, rather than pinpointing the important contributions made along the way to the current shape of a certain result. We stress that this is not a comprehensive survey and the choices of reflect bias. 

\me In the first part we are concerned with theoretical aspects: definitions and relations. In the second part we focus on the practical aspect of computing singularity invariants and we review certain classes of singularities.

\me I would like to thank  A. Dimca, G.-M. Greuel, K. Sugiyama, K. Takeuchi, and W. Veys for their help, comments, and suggestions. Also I would like to thank Universit\'{e} de Nice for their hospitality during the writing of  this article.

\section{Theoretical aspects.}

\subsection{Topology.} 
\

\me
\ni
{\bf Milnor fiber and monodromy.} Let $f$ be a hypersurface singularity germ at the origin in $\bC^n$. 

\medskip
Let
$$M_{t}:=f^{-1}(t)\cap B_\epsilon ,$$
where $B_\epsilon$ is a ball of radius $\epsilon$ around the
origin. Small values of $\eps$ and even smaller values of 
$|t|$  do not change the diffeomorphism class of
$M_t$, the {\it Milnor fiber of $f$ at $0$} \cite{M, Ku}.

\medskip
Fix a Milnor fiber $M_t$ and let $$M_{f,0}:=M_t.$$ The cohomology groups
$H^i(M_{f,0},\bC)$  admit an action $T$
called {\it monodromy} generated by going once around a loop
starting at $t$ around $0$. The eigenvalues of the monodromy action $T$ are roots of
unity, \cite{M, Ku}.

\me
The {\it monodromy zeta function of $f$ at $0$} is
$$
Z^{mon}_{0}(s):=\prod_{j\in\bZ} \det (1-sT, H^j(M_{f,0},\bC))^{(-1)^j}.
$$ 

\me
The {\it $m$-th Lefschetz number of $f$ at $0$} is
$$
\Lambda (T^m):=\sum_{j\in\bZ} (-1)^j \text{Trace}\, (T^m,H^j(M_{f,0},\bC)).
$$
These numbers recover the monodromy zeta function: if $\Lambda (T^m)=\sum_{i|m}s_i$ for $m\ge 1$, then $Z^{mon}_0(s)=\prod_{i\ge 1}(1-t^i)^{s_i/i}$, \cite{Di}.

\me
When $f$ has an isolated singularity, 
$$
\dim_\bC H^j(M_{f,0},\bC)=
\left \{\begin{array}{ll}
0 &\text{ for }j\ne 0, n-1,\\
1 &\text{ for }j=0,\\
\dim_\bC {\bC[[x_1,\ldots ,x_n]]/\left (\frac{\pa f}{\pa x_1},\ldots ,\frac{\pa f}{\pa x_n}\right)} & \text{ for }j=n-1.
\end{array}
\right.
$$
The last value for $j=n-1$ is denoted $\mu(f)$ and called the {\it Milnor number of $f$}, \cite{M, Ku}.

\me
The most recent and complete textbook on the basics, necessary to understand many of the advanced topics here is \cite{GLS}.

\me\ni{\bf Constructible sheaves.}
Let $X$ be a nonsingular complex variety and $Z$ a closed subscheme. 

\me
The Milnor fiber and the monodromy can be generalized to this setting. Let $D^b_c(X)$ be the derived category of bounded complexes of sheaves  of
$\bC$-vector spaces with constructible cohomology in the analytic topology of $X$, \cite{Di}.

\medskip
If $Z$ is a hypersurface given by a regular function $f$, one has
Deligne's {\it nearby cycles functor}.  This is the composition of derived functors
$$\psi_f := i^*p_*p^*  : D^b_c(X) \ra D^b_c(Z_{red}) ,$$
where $i$ is the inclusion of $Z_{red}$ in $X$, $\tilde{\bC}^*$ is the universal cover of $\bC^*$, and $p:X\times_{\bC}\tilde{\bC}^* \ra X$ is the natural projection. If $i_x$ is the inclusion of a point $x$ in $Z_{red}$ and $M_{f,x}$ is the Milnor fiber of $f$ at $x$, then
$$
H^i(i_x^*\psi_f\bC_X)=H^i(M_{f,x},\bC)
$$
and there is an induced action recovering the monodromy, \cite{Di}.

\medskip
When $Z$ is closed subscheme one has {\it Verdier's specialization functor $\Sp_Z$.} This is defined by using $\psi_t$, where $t:\cX\ra\bC$ is the deformation to the normal cone of $Z$ in $X$.
This functor recovers the nearby cycles functor $\psi_f$ in the case when $Z=\{f=0\}$, \cite{Ve}. Another functor,  also recovering the nearby cycles functor, seemingly depending on equations $f=(f_1,\ldots ,f_r)$ for $Z$, is {\it Sabbah's specialization functor $^A\psi_f$}. This is defined by replacing in the definition of the nearby cycles functor $\bC$ and $\bC^*$ with $\bC^r$ and $(\bC^*)^r$, respectively, \cite{Sab}. It would be interesting to understand the differences between these two functors.

\subsection{Analysis.}	
\

\me\ni{\bf Asymptotic expansions.} 
Let $f$ be a hypersurface germ with an isolated singularity in
$\bC^n$.

\medskip
Let $\sigma$ be a top relative holomorphic form on the Milnor fibration $M\ra S$, where $S$ is a small disc and $M=\cup_{t\in S}M_t$. Let $\delta_t$ be a continuous family of  cycles in $H_{n-1}(M_t,\bC)$. Then 
$$
\lim_{t\ra 0}\int_{\delta_t}\sigma = \sum_{\alpha\in\bQ,
k\in\bN} a(\sigma,\delta,\alpha,k)\cdot t^\alpha (\log t)^k$$
where $a$ are constants. The infimum of rational numbers $\al$ that can appear in such expansion for some $\sigma$ and $\delta$ is Arnold's {\it complex oscillation index}. This is an analytic invariant, \cite{AGV, Ko}.

\me\ni{\bf $L^2$-multipliers.}
Let $f$ be a collection of polynomials
$f_1,\ldots ,f_r$ in $\bC[x_1,\ldots ,x_n]$, and let  $c$ be a positive real number.

\medskip
The {\it multiplier ideal of $f$ with coefficient $c$} of Nadel is
the ideal sheaf $\cJ(f^c)$ consisting locally of holomorphic functions $g$ such that $|g|^2/(\sum_i|f_i|^2)^c$ is locally integrable. This is a coherent ideal sheaf. The intuition behind this analytic invariant is: the smaller the
multiplier ideals are, the worse the singularities of the zero locus of $f$ are, \cite{La}. 

\medskip
The smallest $c$ such that $\cJ(f^c)\ne \cO_X$, i.e.
$1\not\in\cJ(f^c)$, is called the {\it log canonical
threshold} of $f$ and is denoted $\lct (f)$.

\me Log canonical thresholds are a special set of numbers: for a fixed $n$, the set $\{\lct(f)\ |\ f\in\bC[x_1,\ldots ,x_n]  \}$ satisfies the ascending chain condition, \cite{DEM}.

\medskip When $f$ is only one polynomial with an isolated singularity, the log canonical threshold 
 coincides with $1 + $ Arnold's complex oscillation index \cite{Ko}. 
 
 \medskip
The definition of the multiplier ideal generalizes and patches up to define, globally on a nonsingular variety $X$ with a
subscheme $Z$, a {\it multiplier ideal sheaf} $\cJ(X,c\cdot Z)$ in $\cO_X$. In fact, the multiplier ideal $\cJ(X,c\cdot Z)$ depends only on the integral closure of the ideal of $Z$ in $X$. One has similarly a {\it log canonical threshold} for $Z$ in $X$, denoted $\lct(X,Z)$, \cite{La}.

\subsection{Geometry.}
\

\me\ni{\bf Resolution of singularities.} Let $X$ be a nonsingular complex variety and $Z$ a closed subscheme. 

\medskip
Let $\mu:Y\ra X$ be a {\it log
resolution} of $(X,Z)$.  This means that $Y$ is nonsingular, $\mu$ is birational and proper, and the inverse image of $Z$ together with the support of the determinant of the Jacobian of $\mu$ is a simple normal crossings divisor. This exists by Hironaka. Denote by $K_{Y/X}=\sum_{i\in S} k_iE_i$ the divisor given by
the determinant of the Jacobian of $\mu$. Denote by $E=\sum_{i\in S} a_iE_i$ the divisor
in $Y$ given by $Z$. Here $E_i$ are irreducible divisors. For $I\subset S$, let $E_I^o:= \cap_{i\in I}E_i-\cup_{i\not\in S}E_i$.

\me
Let $c\in\bR_{>0}$. Then Nadel's multiplier ideal equals 
$$\cJ(X,c\cdot Z)= \mu_*\cO_Y(K_{Y/X}-\rndown{c\cdot E}).$$
Here
$\rndown{\cdot}$ takes the round-down of the coefficients of the
irreducible components of a divisor, \cite{La}. 

\me In particular, the log canonical threshold is given by
\begin{align}\label{eqlct}
\lct (X,Z)=\min_i \left\{\frac{k_i+1}{a_i}\right\}.
\end{align}

\me 
When $Z=\{f=0\}$  is a hypersurface and $x\in Z$ is a point, the monodromy zeta function at $x$ and the Lefschetz numbers can be computed from the log resolution by {\it A'Campo formula}:
$$
\Lambda (T^m)=\sum_{a_i|m} a_i\cdot\chi (E_i^o\cap \mu^{-1}(x)),
$$
where $\chi$ is the topological Euler characteristic, \cite{Di}.

\me
One can imitate the construction via log resolutions to define {\it multiplier ideals for any linear combination of subschemes}, or equivalently, {\it of ideals}:

 $$\cJ(X,c_1\cdot Z_1+\ldots +c_r\cdot Z_r)=\cJ(X,I_{Z_1}^{c_1}\cdot\ldots\cdot I_{Z_r}^{c_r}).$$

\me
If $X=\bC^n$ and the ambient dimension $n$ is $<3$, every integrally closed ideal is a multiplier ideal \cite{LW}. This is not so if $n\ge 3$, \cite{LL}.

\me
The {\it jumping numbers of $Z$ in $X$} are those numbers $c$ such that $$\cJ(c\cdot
Z)\ne\cJ((c-\epsilon)\cdot Z)$$ for all $\eps>0$. The log canonical threshold $\lct(X,Z)$ is the smallest jumping number. The list of jumping numbers is another numerical analytic invariant of the singularities of $Z$ in $X$. The list contains finitely many numbers in any compact interval, all rational numbers, and is periodic. If $\lct(f)=c_1<c_2 <\ldots $ denotes the list of jumping numbers then 
$$c_{i+1}\le c_1+c_i .$$  The standard reference for jumping numbers and multiplier ideals is \cite{La}.

\me
For a point $x$ in $Z$, the {\it inner jumping multiplicity of $c$ at $x$} is the vector space dimension
$$
m_{c,x}:=\dim_\bC \cJ(X,(c-\eps)\cdot Z) / \cJ(X,(c-\eps)\cdot Z+\delta\cdot \{x\}),
$$
where $0<\eps\ll\delta\ll 1$. This multiplicity measures the contribution of the singular point $x$ to the jumping number $c$, \cite{B-HS}.

\me Another interesting singularity invariant is the  Denef-Loeser {\it topological zeta function}. This is the rational function of complex variable $s$ defined as
$$
Z^{top}_Z (s):= \sum_{I\subset S}\chi (E_I^o) \cdot \prod_{i\in I}\frac{1}{a_is+k_i+1}.
$$
This is independent of the choice of log resolution, \cite{DL}. In spite of the name, $Z^{top}_f(s)$ is not a topological invariant, \cite{AB}.

\me When $Z=\{f=0\}$ is a hypersurface, the {\it Monodromy Conjecture} states that if $c$ is pole of the topological zeta function, then $e^{2\pi i c}$ is an eigenvalue of the Milnor monodromy of $f$ at some point in $f^{-1}(0)$, \cite{DL}. A similar conjecture, using Verdier's specialization functor $\Sp_Z$, can be made when $Z$ is not a hypersurface, \cite{PV}.

\me\ni{\bf Mixed Hodge structures.} 
Let $X$ be a nonsingular complex variety and $Z$ a closed subscheme. 

\me
The topological package consisting of the Milnor fibers, monodromy, nearby cycles functor, specialization functor can be enhanced to
take into account natural mixed Hodge structures, \cite{S}. 

\me
Consider the case when $Z$ is a hypersurface given by one polynomial $f\in \bC[x_1,\ldots ,x_n]$ with the origin included in the singular locus. The {\it Hodge spectrum of $f$
at $0$} of Steenbrink is
\begin{equation*}
\Sp(f,0)=\sum_{c>0}n_{c,0}(f)\cdot t^c,
\end{equation*}
where  the {\it spectrum multiplicities}
$$n_{c,0}(f):=\sum_{i\in\bZ}(-1)^{n-1-i}\dim_\bC \Gr _F ^{\rndown{n-c}}\wti{H}^{i}(M_{f,0},\bC)_{e^{-2\pi ic}}$$ 
record the generalized Euler
characteristic on the $\rndown{n-c}$-graded piece of the Hodge
filtration on the $\exp (-2\pi ic)$-monodromy eigenspace on the reduced cohomology of the Milnor fiber. These invariants can be refined by considering the weight filtration as well, \cite{Ku}. 

\me
In the case of isolated hypersurface singularities, the spectrum  recovers the Milnor number
$$
\mu(f)=\sum_c n_{c,0}(f)
$$
and, by M. Saito, the {\it geometric genus} of the singularity 
$$
\sum_{0<c\le 1}n_{c,0}(f)= p_g(f) :=\left\{
\begin{array}{ll}
\dim_\bC(R^{n-2}p_*\cO_{\ti{Z}})_0 & \text{ if }n\ge 3 ,\\
\dim_\bC (p_*\cO_{\ti{Z}} / \cO_{Z})_0 & \text{ if }n=2,
\end{array}\right .
$$
where $p:\ti{Z}\ra Z$ is a log resolution of $Z$, \cite{Ku}. The spectrum also satisfies a symmetry $n_{c,0}(f)=n_{n-c,0}(f)$, and a semicontinuity property. It is thus useful in the classification of such singularities, \cite{Ku}. The smallest spectral number $c$  equals the log canonical threshold $\lct (f)$. Let $c_1\le \ldots\le c_{\mu(f)}$ denote the list of spectral numbers counted with the spectrum multiplicities. An open question is {\it Hertling's Conjecture} stating that 
$$
\frac{1}{\mu(f)}\sum_{i=1}^{\mu(f)}\left (   c_i - \frac{n}{2}\right )^2 \le \frac{c_{\mu(f)}-c_1}{12} .
$$
This has been solved for quasi-homogeneous singularities \cite{He}, where equality holds. This is due to a duality with the spectrum of the Milnor fiber at infity, for which in general a similar conjecture is made but with reversed sign, \cite{Di-hert}. Other solved cases are:  irreducible plane curves \cite{S-ex} and Newton nondegenerate polynomials of two variables \cite{Bre}. Another open question is {\it Durfee's Conjecture} of \cite{Dur} that for $n=3$,
$$
6p_g(f)\le \mu(f) .
$$
This was shown to be true in the following cases: quasi-homogeneous \cite{XY}, weakly elliptic, $f=g(x,y)+z^N$ \cite{As, Nem}, double point \cite{Tom}, triple point \cite{As-tri}, absolutely isolated \cite{MH-abs}.

\me
The jumping numbers are also related to Milnor fibers and monodromy. If the singularity is isolated, the spectrum recovers all the jumping numbers in $(0,1)$. In general, when the singularities are not necessarily isolated, we have more precisely that the spectrum multiplicities for $c\in (0,1]$ are computed in terms of the {inner jumping multiplicities} of jumping
numbers: $m_{c,x}(f)=n_{c,x}(f)$,  \cite{B-HS}.

\me A sufficient condition for symmetry of the spectrum of a homogeneous polynomial in the non-isolated case is given in \cite{Di-Mon}-Prop. 4.1.

\me
The Hodge spectrum has a generalization to any subscheme $Z$ in a nonsingular variety $X$, using Verdier's specialization functor and M. Saito's mixed Hodge modules. There is a relation between the multiplier ideals and the specialization functor, \cite{DMS}.

\me\ni{\bf Jets and arcs.} 
Let $X$ be a nonsingular complex variety of dimension $n$ and $Z$ a closed subscheme. 

\me
The {\it scheme of $m$-jets} and the {\it arc space} of $Z$ are
\begin{align*}
Z_m &:=\Hom ({\rm{Spec}}\; \bC[t]/(t^{m+1}) , Z)\quad \text{respectively}\\
Z_\infty &:=\Hom ({\rm{Spec}}\; \bC[[t]], Z).
\end{align*}

\me
Jets compute log canonical thresholds by Musta\c{t}\u{a}'s formula \cite{Mu}:
$$
\lct (X,Z)=\min _m\left\{\frac{\codim (Z_m,X_m)}{m+1}\right\}.
$$

\me If $Z=\{f=0\}$ is a hypersurface, one has  the {\it Denef-Loeser motivic zeta function}:
$$
Z^{mot}_f(s):=\sum_{m\ge 1}[X_{m,1}][\bA^1]^{-mn}s^m,
$$
where $[.]$ denotes the class of a variety in an appropriate Grothendieck ring, and $X_{m,1}$ consists of the $m$-jets $\phi$ of $X$ such that $f(\phi)=t^m$.  The motivic zeta function is a rational function. The monodromy zeta function, the Hodge spectrum, and the topological zeta function can be recovered from the motivic zeta function. Thus these singularity invariants can be computed from jets. In fact, one has the {\it motivic Milnor fiber} $$\cS_f:=-\lim_{s\ra \infty}Z^{mot}_f(s),$$ which is a common generalization of the monodromy zeta function and of the Hodge spectrum, \cite{DL}. 

\me
The motivic zeta function can be defined also when $Z$ is a closed subscheme with equations $f=(f_1,\ldots ,f_r)$. The {\it Monodromy Conjecture} can be stated for the motivic zeta function and implies the previous version, \cite{Ni}. The analog of the motivic Milnor fiber $\cS_f$  is related  in this case with Sabbah's specialization functor $^A\psi_f$: it recovers the generalization via $^A\psi_f$ of the monodromy zeta function, \cite{Gu}.

\me
The biggest pole of $Z^{mot}_f(s)$ gives the negative of the log canonical threshold, \cite{HN}-p.18.

\me
The motivic zeta function is a {\it motivic integral}. Without explaining what this is, a motivic integral enjoys a change of variables formula. In practice this means that a motivic integral can be computed from a log resolution. From this very advanced point of view, one can see more naturally A'Campo's formula for the monodromy zeta function and the connection between jumping numbers and the Hodge spectrum, \cite{DL}.

\me
The change of variables formula can be streamlined and the motivic integration eliminated. The {\it $m$-th contact locus of $Z$ in $X$}  is the subset of $X_\infty$ consisting of arcs of order $m$ along $Z$. The contact loci can also be expressed in terms of log resolutions and exceptional divisors in log resolutions give rise to components of contact loci. In many cases this also gives a back-and-forth pass between arc-theoretic invariants and invariants defined by log resolutions, such as Musta\c{t}\u{a}'s result on the log canonical threshold. Multiplier ideals and jumping numbers can also be interpreted arc-theoretically, \cite{ELM}.

\me 
If $Z$ is a normal local complete intersection variety, then $Z_m$ is equidimensional (respectively irreducible, normal) for every $m$ if and only if $Z$ has log canonical (canonical, terminal) singularities, \cite{EM}.

\subsection{Algebra.}
\

\me\ni{\bf Riemann-Hilbert correspondence.} Let $X$ be nonsingular complex variety of dimension $n$.

\me The {\it sheaf of algebraic differential operators} $\cD_X$ is locally given in affine coordinates by the Weil algebra $\bC[x_1,\ldots ,x_n,\pa /\pa x_1 ,\ldots ,\pa/\pa x_n]$. An important class of (left) $\cD_X$-modules consists of those {\it regular and holonomic}. Let $D^b_{rh}(\cD_X)$ be the bounded derived category of complexes of $\cD_X$-modules with regular holonomic cohomology, \cite{Bo}.

\me One of the main reasons why the theory of $\cD$-modules has become important recently is because of its suitability for computer calculations.

\me The topological package, consisting of the bounded derived category of constructible sheaves $D^b_c(X)$ and the natural functors attached to it, has an algebraic counterpart. There is a well-defined functor
$$
DR: D^b_{rh}(\cD_X) \ra D^b_c(X)
$$
which is an equivalence of categories commuting with the usual functors,  \cite{Bo}. 

\me The $\cD$-module theoretic counterpart of the nearby cycles functor $\psi_f$, hence of the Milnor monodromy of $f$, is achieved by the {\it $V$-filtration along $f$} of Malgrange-Kashiwara. For $c\in (0,1)$,
$$
\psi_{f,\lam} \bC_X[-1]=DR (\Gr _V^c\wti{\cO_X}),
$$
where $\lam=e^{-2\pi ic}$, $\psi_f=\oplus_{\lam'} \psi_{f,\lam '}$ is the functor decomposition corresponding to the eigenspace decomposition of the semisimple part of the Milnor monodromy, and $\wti{\cO_X}=\cO_X[\pa_t]$ is the $\cD$-module push-forward of $\cO_X$ under the graph embedding of $f$, \cite{B-V}.

\me
In algebraic geometry, integral (co)homology groups are endowed with additional structure: mixed Hodge structures, \cite{PS}. The modern point of view is M. Saito's theory of {\it mixed Hodge modules}. The derived category of mixed Hodge modules $D^b(MHM(X))$ has natural forgetful functors to $D^b_{rh}(\cD_X)$ and $D^b_c(X)$ and recovers Deligne's mixed Hodge structures on the usual (co)homology groups.  When $Z$ is a closed subvariety of $X$, the Verdier specialization functor $Sp_Z$ also exists in the framework of mixed Hodge modules,  \cite{S}.

\me 
Let $Z$ be a closed subscheme of $X$, and let $\wti{\cO_X}$ be the $\cD$-module push-forward of $\cO_X$ under the graph embedding of a set of local generators of the ideal of $Z$ in $X$. The smallest nontrivial piece of the Hodge filtration of the $V$-filtration on $\wti{\cO_X}$ gives the multiplier ideals:
$$
 \cJ(X,(c-\eps)\cdot Z)=F^{n-1}V^c\wti{O_X},
$$
with $0<\eps\ll 1$. This is another point of view on the relation between multiplier ideals, mixed Hodge structures, and Milnor monodromy, \cite{B-Comp}.

\me\ni{\bf $b$-functions.} Let $X$ be a nonsingular complex variety of dimension $n$ and $Z$ a closed subscheme.

\me
Suppose first that $Z$ is given by an ideal  $f=(f_1,\ldots , f_r)$  with $f_i\in\bC[x_1,\ldots ,x_n]$. Let $g$ be another polynomial in $n$ variables. The {\it generalized $b$-function of $f$ twisted by $g$}, also called the {\it generalized Bernstein-Sato polynomial of $f$ twisted by $g$} and denoted $b_{f,g}(s)$, is the nonzero monic polynomial of minimal degree among those $b\in\bC[s]$ such that
$$
b(s_1+\ldots + s_r)g\prod_{i=1}^r f_i^{s_i}  = \sum_{k=1}^r P_k(gf_k\prod_{i=1}^r f_i^{s_i}),
$$
for some algebraic operators $P_k\in\bC[x_1,\ldots ,x_n, \pa/\pa x_1,\ldots ,\pa/\pa x_n][s_{ij}]_{1\le i,j\le r}$, where $s_{ij}$ are defined as follows. First, let the operator $t_i$ act by leaving $s_j$ alone if $i\ne j$, and replacing $s_i$ with $s_i+1$. For example: $t_j\prod_{i=1}^r f_i^{s_i}=f_j\prod_{i=1}^r f_i^{s_i}$. Then $s_{ij}:=s_it_i^{-1}t_j$. The generalized $b$-function is independent of the choice of local generators $f_1,\ldots ,f_r$ for the ideal of $Z$, \cite{B-Comp}.

\me
The {\it  $b$-function of the ideal $f$} is  $b_f(s):=b_{f,1}(s)$. When $Z$ is a hypersurface, $b_f(s)$ is the usual $b$-function of Bernstein and Sato, satisfying the relation
$$
b_f(s)f^s=Pf^{s+1}
$$
for some operator $P_k\in\bC[x_1,\ldots ,x_n, \pa/\pa x_1,\ldots ,\pa/\pa x_n][s]$.

\me For a scheme $Z$, the {\it $b$-function of $Z$} is the polynomial $b_{Z}(s)$ obtained by replacing $s$ with $s-\codim (Z,X)$ in the lowest common multiple of the polynomials $b_f(s)$ obtained by varying local charts of a closed embedding of $Z$ into a nonsingular $X$. This polynomial depends only on $Z$. The roots of $b_{f,g}(s)$ and  $b_f(s)$ are negative rational numbers,
 \cite{B-Comp}.

\me
For simplicity, say $X$ is affine from now and the ideal of $Z$ is $f$.  The $b$-function recovers the monodromy eigenvalues, as observed originally by Malgrange and Kashiwara. If $Z$ is a hypersurface, the set consisting of $e^{2\pi i c}$, where $c$ are roots of $b_f(s)$, is the set of eigenvalues of the Milnor monodromy at points along $Z$. In higher codimension, a similar statement holds for eigenvalues related with the specialization functor $\Sp _Z$, \cite{B-Comp}.

\me The {\it Strong Monodromy Conjecture} states that if $c$ is a pole of the topological zeta function $Z^{top}_f(s)$ of the ideal $f$, then $b_f(c)=0$. It can be stated for the motivic zeta function as well. It implies the Monodromy Conjecture.

\me The biggest root of the $b$-function $b_f(s)$ of the ideal of $Z$ is the negative of the log canonical threshold of $(X,Z)$, \cite{Ko, B-Comp}.

\me
If $\lct (X,Z)\le c< \lct (X,Z)+1$ and $c$ is a jumping number of $Z$ in $X$, then $b_f(-c)=0$, \cite{ELSV, B-Comp}.

\me Generalized $b$-functions recover multiplier ideals, \cite{B-Comp}:
$$
\cJ(X,c\cdot Z) =_{loc} \{ g\in\cO_X\ |\ c<\al \text{ if } b_{f,g}(-\al)=0  \}.
$$

\me
Generalized $b$-functions are related to the $V$-filtration along $f$. More precisely, $b_{f,g}(s)$ is the minimal polynomial of the action of 
$$s=-(\pa_1t_1+\ldots \pa _rt_r)\text{ on }V^0\cD_{Y}(g\otimes 1)/V^1\cD_{Y}(g\otimes 1).$$
Here $Y=X\times \bC^r$, the coordinate functions on $\bC^r$ are $t_1,\ldots ,t_r$, the operator $\pa_j$ is $\pa/\pa t_j$, $\wti{\cO_X}=\cO_X[\pa_1,\ldots ,\pa_r]$ is viewed as a $\cD_Y$-module via the graph embedding of $f$, $g\otimes 1\in\wti{\cO_X}$, and $V^i\cD_Y$ consists of operators $P$ in $\cD_Y$ such that $P(t_1,\ldots ,t_r)^i\subset (t_1,\ldots ,t_r)^{i+j}$ for all $j\in\bZ$, \cite{B-Comp}.

\me 
The $b$-function of a polynomial can sometimes be calculated via {\it microlocal calculus}. This method has been successful for computation of relative invariants of irreducible regular prehomogeneous vector spaces, see below, \cite{Ka-book}.

\subsection{Arithmetic.}
\

\me\ni{\bf $K$-log canonical thresholds.} Let $f\in K[x_1,\ldots ,x_n]$ be a polynomial with coefficients in a  complete field $K$ of characteristic zero. 

\me
By Hironaka, $f$ admits a $K$-analytic log resolution of the zero locus of $f$ in $K^n$. In general this might be too small to be a log resolution over an algebraic closure of $K$. Imitating (\ref{eqlct}), one has the {\it $K$-log canonical threshold of $f$}, which we denote $lct_K(f)$. We have that $\lct (f)=\lct_\bC (f)$ for any embedding $K\subset\bC$, but in general $\lct(f)\le \lct_K(f)$. This can be generalized to the case of an ideal $f$ of polynomials with coefficients in $K$, \cite{VZ}.

\me If $K=\bR$ one can also define the {\it real jumping numbers} of $f$. It can happen that a real jumping number is not an usual jumping number. However, any real jumping number smaller than $lct_\bR\, (f)+1$ is a root of $b_f(-s)$, \cite{S-real}.

\me\ni{\bf $p$-adic local zeta functions.} 
Let $p$ be a prime number and $K$ a finite extension field of $\bQ_p$ with a fixed embedding into $\bC$. Let $f$ be an ideal of polynomials in $K[x_1,\ldots ,x_n]$.

\me
If $f$ is a single polynomial  with coefficients in $\bQ_p$, the {\it Igusa $p$-adic local zeta function of $f$} is defined for a character $\chi$ on the units of $\bZ_p$ as
$$ Z^{p}_{f,\chi} (s):=\int_{\bZ^n_p} |f(x)|^s\chi(ac(f(x)))dx\ ,$$
where $|t|=p^{-ord _pt}$, $ac(t)=|t|t$, and $dx$ is the Haar measure normalized such that the measure of $p^{m_1}\bZ_p\times\ldots\times p^{m_n}\bZ_p$ is $p^{-(m_1+\ldots +m_n)}$. Let $Z^p_f(s)=Z^p_{f,1}(s)$, \cite{De}.

\me
The poles of $Z^p_{f,\chi}(s)$ determine by \cite{Ig} the asymptotic expansion as $|t|\ra 0$ of the numbers
$$
N_m(t):=\{ x\in (\bZ/p^m\bZ)^n\ |\ f(x)\equiv t\text{ mod }p^m   \}\quad (m\gg 0).
$$

\me
The definition of $Z^p_f(s)$ can be made more generally for an ideal $f$ of polynomials with coefficients in a finite extension $K$ of $\bQ_p$. The $p$-adic local zeta functions are rational, \cite{HMY}.

\me If $K=\bQ_p$, the $K$-log canonical threshold is determined by the numbers $N_m:=N_m(0)$:
$$ \lim_{m\ra \infty}
(N_m)^{1/m} = p^{\;n-\lct_K(f)} .$$
A similar statement holds for a finite extension $K$ of $\bQ_p$, see \cite{VZ}.

\me
If $c$ is the pole of $Z^p_f(s)$ with the biggest real part, then $\lct_K(f)=-Re(c)$, \cite{VZ}.

\me
The motivic zeta function $Z^{mot}_f(s)$ of Denef-Loeser determines the $p$-adic local zeta function $Z^p_{f}(s)$, \cite{DL}. 

\me
If $f$ is a single polynomial, Igusa's original {\it Monodromy Conjecture} states that if $c$ is a pole of $Z^p_f(s)$ then $e^{2\pi iRe(c)}$ is an eigenvalue of the Milnor monodromy of $f_\bC$ at some point of $f_\bC^{-1}(0)$, where $f_\bC$ is $f$ viewed as a polynomial with complex coefficients. If $f$ is an ideal defining a subscheme $Z$, the Monodromy Conjecture is stated via Verdier's specialization functor $\Sp _Z$. The {\it Strong Monodromy Conjecture} states that if $f$ is an ideal and $c$ is a pole of $Z^p_f(s)$, then $b_f(Re(c))=0$, \cite{Ig-book}.

\me\ni{\bf Test ideals.}
Let $p$ be a prime number and $f$ be an ideal of polynomials in $\bF_p[x_1,\ldots ,x_n]$.

\me
The Hara-Yoshida {\it test ideal of $f$ with coefficient $c$}, where $c$ is positive real number, is
$$
\tau (f^c):=\left(  f^{\rndup{cp^e}}  \right )^{[1/p^e]},\quad e\gg 0,
$$
where for an ideal $I$, the ideal $I^{[1/p^e]}$ is defined as follows. This is the unique smallest ideal $J$ such that $I\subset \{ u^{p^e}\ |\ u\in J\}$, \cite{BMS}. 

\me
The {\it $F$-jumping numbers} of $f$ are the positive real numbers $c$ such that 
$$
\tau(f^c)\ne \tau (f^{c-\eps})
$$
for all $\eps>0$. The Takagi-Watanabe {\it $F$-pure threshold} of $f$ is the smallest $F$-jumping number and is denoted $\fpt(f)$, \cite{BMS}.

\me
Test ideals, $F$-jumping numbers, and the $F$-threshold are positive characteristic analogs of multiplier ideals, jumping numbers, and respectively, the log canonical threshold, \cite{HY}. More precisely, let now $f$ be an ideal of polynomials in $\bQ[x_1,\ldots ,x_n]$. For large prime numbers $p$, let $f_p\subset \bF_p[x_1,\ldots ,x_n]$ denote the reduction modulo $p$ of $f$.  Fix $c>0$. Then for $p\gg 0$,
$$
\tau (f_p^c)=\cJ(f^c)_p
$$
and
$$
\lim_{p\ra\infty} \fpt (f_p)=\lct (f).
$$

\me
The {\it Hara-Watanabe Conjecture} \cite{HW} states that there are infinitely many prime numbers $p$ such that for all $c>0$, 
$$\tau(f^c_p)=\cJ(f^c)_p.$$

\me The list of $F$-jumping numbers enjoys similar properties as the list of jumping numbers: rationality, discreteness, and periodicity, \cite{BMS}. However, in any ambient dimension, every ideal is a test ideal, in contrast with the speciality of the multiplier ideals, \cite{MY}.

\me
There are results connecting test ideals with $b$-functions. If $f\in\bQ[x_1,\ldots,x_n]$ is a single polynomial and $c$ is an $F$-jumping number of the reduction $f_p$ for some $p\gg 0$, then $\rndup{cp^e}-1$ is a root of $b_f(s)$ modulo $p$, \cite{MTW, M-b}.

\subsection{Remarks and questions.}
\

\me\ni{\bf Answer to the trivia question.} In how many ways can the log canonical threshold of a polynomial be computed? We summarize some of the things we have talked about so far. The lct can be computed, theoretically, via: the $L^2$ condition, the orders of vanishing on a log resolution, the growth of the codimension of jet schemes, the poles of the motivic zeta function, the $b$-function, and the test ideals. If a log resolution over $\bC$ is practically the same as a $K$-analytic log resolution over a $p$-adic field $K$ containing all the coefficients of the polynomial, i.e. if $\lct_K(f)=\lct(f)$, then there are two more ways: via the poles of $p$-adic local zeta functions and via the asymptotics of the number of solutions modulo $p^m$. If the singularity is isolated it can also be done via Arnold's complex oscillation index and via the Hodge spectrum. So, 6+2+2 ways.

\me\ni{\bf What topics were left out.} May topics are left out from this survey: singularities of varieties inside singular ambient spaces, the Milnor fiber at infinity, the characteristics classes point of view on singularities, other invariants such as polar and Le numbers, the theory of Brieskorn lattices, local systems, archimedean local zeta functions, deformations, equisingularity, etc.

\me\ni{\bf Questions.}  We have already mentioned the Monodromy Conjecture and its Strong version, the Hertling Conjecture, the Durfee Conjecture, and the Hara-Watanabe Conjecture.

\me
It is not known how to relate $b$-functions with jets and arcs. In principle, this would help with the Strong Monodromy Conjecture.

 \me 
We know little about the most natural singularity invariant, the {\it multiplicity}. {\it Zariski conjecture} states that if two reduced hypersurface singularity germs in $\bC^n$ are embedded-topologically equivalent then their multiplicities are the same. Even the isolated singularity case is not known, \cite{Ey}. It known to be true for semi-quasihomogeneous singularities: \cite{Gre}, and slightly weaker, \cite{OSh}.

\me We can raise the same question for log canonical thresholds. Can one find an example of two reduced hypersurface singularity germs in $\bC^n$ that are embedded-topologically equivalent but have different log canonical thresholds?

\me Is the biggest pole of the topological zeta function $Z^{top}_f(s)$ of a polynomial $f$ equal to $-\lct (f)$? This true for 2 variables, \cite{V-top}.

\me For a polynomial $f$ with coefficients in a complete field $K$ of characteristic zero, define $K$-jumping numbers and prove the ones $<\lct_K(f)+1$ are roots of $b_f(-s)$, as in the cases when $K$ is $\bC$ or $\bR$.

\me Can microlocal calculus, which provides a method for computation of the $b$-function of a polynomial, be made to work for $b$-functions of ideals ?

\me Let $f=(f_1,\ldots ,f_r)$ be a collection of polynomials. What are the differences between Verdier's and Sabbah's specialization functor for $f$? Does the motivic object $\cS_f$, the higher-codimensional analog of the motivic Milnor fiber of a hypersurface, recover the generalized Hodge spectrum of $f$? 

\me Can the geometric genus $p_g$ of a normal isolated singularity can be recovered from the generalized Hodge spectrum, in analogy with the isolated hypersurface case? This would be relevant to the original, more general form of Durfee's Conjecture, which was stated for isolated complete intersection singularities.

\section{Practical aspects.}

\subsection{General rules.}  We mention some rules that apply for calculation of singularity invariants or help approximate singularity invariants. Whenever a geometric construction is available, one can look for the formula describing the change in a singularity invariant. We have already talked about log resolutions and jet schemes.

\me
An {\it additive Thom-Sebastiani rule} describes a singularity invariant for $f(x)+g(y)$ in terms of the invariants for $f$ and $g$, when $f(x)$ and $g(y)$ are polynomials in two disjoint sets of variables. This rule is  available: for the motivic Milnor fiber, and hence for the monodromy zeta function and the Hodge spectrum, \cite{DL}; for the poles of the $p$-adic zeta functions, \cite{DV}; and for the $b$-function when both polynomials have isolated singularities and $g$ is also quasihomogeneous, \cite{Ya}. 

\me
An {\it additive Thom-Sebastiani rule for ideals} describes a singularity invariant for a sum of two ideals in two disjoint sets of variables. Equivalently, this rule describes a singularity invariant of a product of schemes. This rule is the easiest one to obtain. It is available for example for motivic zeta functions \cite{DL}, multiplier ideals, jumping numbers \cite{La}, $b$-functions \cite{B-Comp}, and test ideals \cite{Ta}.

\me 
A {\it multiplicative Thom-Sebastiani rule} describes a singularity invariant for $f(x)\cdot g(y)$ in terms of the invariants for $f$ and $g$, when $f(x)$ and $g(y)$ are polynomials in two disjoint sets of variables. This rule is available for the Milnor monodromy  of homogeneous polynomials \cite{Di-Mon}-Thm. 1.4, and, in the even greater generality when $f$ and $g$ are ideals, for multiplier ideals and jumping numbers \cite{La}.

\me 
A more general idea is to describe singularity invariants of $F(f_1,\ldots ,f_r)$, where $F$ is a nice polynomial and $f_1,\ldots ,f_r$ are polynomials in distinct sets of variables. For results in this direction for the motivic Milnor fiber see \cite{GLMc,GLMb, GLMa}.

\me
It is hard to say what a {\it summation rule} should be in general. This is available for multiplier ideals \cite{Mu-sum} and test ideals \cite{Ta}:
\begin{align*}
\cJ((f+g)^c)=\sum_{\lam+\mu=c} \cJ(f^\lam\cdot g^\mu),
\end{align*}
and similarly for test ideals, where $f, g$ are ideals of polynomials in the same set of variables. One can ask if a similar rule exists for the Verdier specialization functor or motivic zeta functions.

\me 
A {\it restriction rule} says that an invariant of a hyperplane section of a singularity germ is the same or worse,  reflecting more complicated singularities, than the one of the original singularity. For example, log canonical and $F$-pure thresholds  get smaller upon restriction. Also multiplier ideals \cite{La} and test ideals \cite{HY} get smaller upon restriction. These invariants also satisfy a {\it generic restriction rule} saying that they remain the same upon restriction to a general hyperplane section. This is related to the {\it semicontinuity rule} stating that singularities get worse at special points in a family. The Hodge spectrum of an isolated hypersurface singularity satisfies a semicontinuity property, \cite{Ku}.

\me
For more geometric transformation rules for multiplier ideals see \cite{La}, for test ideals see \cite{BkS, KT}, for nearby cycles functors, motivic Milnor fibers see \cite{Di, GLMb, GLMa}, for jet schemes see \cite{EM-f}, for log canonical thresholds see \cite{dFEM}.

\subsection {Ambient dimension two.}  For a germ of a reduced and irreducible curve $f$ in $(\bC^2,0)$ one has a set of {\it Puiseux pairs} $(k_1,n_1;\ldots ;k_g,n_g)$ defined via a parametrization of the curve

\begin{align*}
y= & \sum_{1\le i\le\lfloor \frac{k_1}{n_1} \rfloor } c_{0,i}x^i\  + \sum_{0\le i\le\lfloor \frac{k_2}{n_2} \rfloor } c_{1,i}x^{(k_1+i)/n_1}\  +\\
& +\ \sum_{0\le i\le\lfloor \frac{k_3}{n_3} \rfloor } c_{2,i}x^{ k_1/n_1+  (k_2+i)/n_1n_2}\ +\ldots \\
&  \ldots +\  \sum_{0\le i } c_{g,i}x^{k_1/n_1+   k_2/n_1n_2+\ldots +    (k_g+i)/n_1\ldots n_g} ,
\end{align*}
where $c_{j,i}\in\bC$, $c_{j,0}\ne 0$ for $j\ne 0$, $k_j, n_j\in\bZ_+$, $(k_j,n_j)=1$, $n_j>1$, and $k_1>n_1$.
The Puiseux pairs determine the embedded topological type. For every plane curve there is a {\it minimal log resolution}, \cite{C}. 

\me
The Hodge spectrum can be written in terms of the Puiseux pairs  for irreducible curves \cite{S-ex}, and in terms of the graph and the vanishing orders of the minimal log resolution for any curves, \cite{LS}.

\me\ni{\it Example.} If $f$ is a irreducible curve germ, the numbers $c<1$ appearing in the Hodge spectrum of $f$, counted with their spectrum multiplicity $n_{c,0}(f)$, are
$$
\left\{\  \frac{1}{n_{s+1}\ldots n_g}\cdot \left(   \frac{i}{n_s}+\frac{j}{w_s} \right) + \frac{r}{n_{s+1}\ldots n_g} \  \right\}
$$
where: $w_1=k_1$, $w_i=w_{i-1}n_{i-1}n_i+k_i$ for $i>1$, $0<i<n_s$, $0<j<w_s$, $0\le r<n_{s+1}\ldots n_g$, and $1\le s\le g$ such that $i/n_s + j/w_s <1$.

\me
Jumping numbers of reduced curves can be written in terms of the graph and the vanishing orders of the minimal log resolution by \cite{B-HS}, which reduced the problem to the Hodge spectrum. The above example also gives the jumping numbers $c<1$ of an irreducible germ together with their inner jumping multiplicities $m_{c,0}(f)$.

\me For any plane curve germ there exists a local system of coordinates such that the log canonical threshold is $1/t$ where $(t,t)$ is the intersection of the boundary of the Newton polytope (see \ref{subNondeg}) with the diagonal line, \cite{Ale, AN}.

\me
Jumping numbers  for a complete ideal of finite colength in two variables are computed combinatorially, \cite{HJ}.

\me
The $b$-function of almost all irreducible and reduced plane curves with fixed Puiseux pairs is 
conjectured to be determined by a precise formula depending on the Puiseux pairs \cite{Y}. There are only partial results on determination of the $b$-function of plane curves, \cite{De, BMT}.

\me The poles of motivic ($p$-adic, topological) zeta functions for any ideal of polynomials in two variables are determined in terms of the graph and the vanishing orders of the minimal log resolution,  \cite{Gu, PV-p, V-top, De}.

\me The Strong Monodromy Conjecture  is proven by Loeser for reduced plane curves, \cite{Lo-2}. The Monodromy Conjecture is proven for ideals in two variables by Van Proeyen-Veys, \cite{PV}.

\me Jet schemes of plane curves are considered in \cite{Mou}.

\subsection {Nondegenerate polynomials.}\label{subNondeg}  The monomials appearing in a polynomial $f$ in $n$ variables determine a set of points in $\bZ_{\ge 0}^n$ whose convex hull is called the {\it Newton polytope} of $f$. The definition of {\it nondegenerate polynomial} is a condition involving the Newton polytope, which can differ in the literature. This is a  condition that expresses in a precise way the fact that the polynomial is general and that is has an explicit log resolution. The following hold under nondegeneracy assumptions.

\me
There are formulas in terms of the Newton polytope for: the Hodge spectrum \cite{St77, Sa} and the monodromy Jordan normal form \cite{ET} when $f$ has isolated singularities ; multiplier ideals and the jumping numbers in general \cite{La}; 

\me\ni{\it Example.} Let $f$ be nondegenerate in the following sense: the  form $df_\sigma$ is nonzero on $(\bC^*)^n\subset\bC^n$, for every face $\sigma$ of the Newton polytope, where $f_\sigma$ is the polynomial composed of the terms of $f$ which lie in $\sigma$. Then Howald showed that for $c<1$ the multiplier ideals $\cJ(f^c)$ are the same as the multiplier ideals $\cJ(I_f^c)$, where $I_f$ is the ideal generated by the terms of $f$. See next subsection for monomial ideals.

\me 
The poles of $p$-adic zeta functions are among a list determined explicitly by the Newton polytope, \cite{DH, Z}; the same holds for nondegenerate maps $f=(f_1,\ldots , f_r)$, \cite{VZ}. The motivic zeta function and the motivic Milnor fiber are considered in \cite{Gu}.

\me
The Strong Monodromy Conjecture is proved for nondegenerate polynomials satisfying an additional condition, by displaying certain roots of the $b$-function, \cite{Lo-nondeg}.

\subsection{Monomial ideals.} A {\it monomial ideal} is an ideal of polynomials generated by monomials.
The {\it semigroup} of an ideal $I\subset\bC[x_1,\ldots ,x_n]$ is the set $\{ u \ |\ x^u\in I\}$. The convex hull of this set is the {\it Newton polytope} $P(I)$ of the ideal. The Newton polytope of a monomial ideal equals the one of the integral closure of the ideal.

\me
The Newton polytope of a monomial ideal determines explicitly the Hodge spectrum \cite{DMS}, the multiplier ideals and the jumping numbers, \cite{La}; the test ideals and $F$-jumping numbers, \cite{HY}; and the $p$-adic zeta function, \cite{HMY}.

\me\ni{\it Example.} Howald's formula for the multiplier ideals of a monomial ideal $I$ is
$$
\cJ(I^c)=\langle x^u\ |\ u+\text{\bf 1}\ \in\ \text{Interior} (cP(I))  \rangle.
$$

\me
The  $b$-function of a monomial ideal has been computed in terms of the semigroup of the ideal. In general, the $b$-function cannot be determined by the Newton polytope alone, \cite{B-mon1,B-mon2}.

\me 
The Strong Monodromy Conjecture  is checked for monomial ideals, \cite{HMY}.

\me The geometry of the jet schemes of monomial ideals is described in \cite{GS, Yu}.

\subsection {Hyperplane arrangements.} Let $K$ be a field. A {\it hyperplane arrangement} $D$ in $K^n$ is a possibly nonreduced union of hyperplanes of $K^n$. An invariant of $D$ is {\it combinatorial} if it only depends on the lattice of intersections of the hyperplanes of $D$ together with their codimensions.  

\me Blowing up the intersections of hyperplanes gives an explicit log resolution. There is also a minimal resolution, \cite{DP}.

\me Jet schemes of hyperplane arrangements are considered in \cite{Mu-arr}. Multiplier ideals are also considered here, see also \cite{Tei}.

\me A current major open problem in the theory of hyperplane arrangements is the combinatorial invariance of the Betti numbers of the cohomology of Milnor fiber, or stronger, of the dimensions of the Hodge pieces. The simplest unknown case is the cone over a planar line arrangement with at most triple points \cite{CoLi}.

\me
The jumping numbers and the Hodge spectrum of a hyperplane arrangement are explicitly determined combinatorial invariants, \cite{B-SpAr}.

\me\ni{\it Example.} Let $f\in\bC[x,y,z]$ be a homogeneous reduced product of $d$ linear forms. This plane arrangement is a cone over a line arrangement $D\in\bP^2$. The Hodge spectrum multiplicities are:
$$
n_{c,0}(f)=0, \text{ if } cd\not\in\bZ;
$$
$$
n_{c,0}(f)=\binom{i-1}{2}-\sum_{m\ge 3}\nu_m\binom{\lceil \frac{im}{d} \rceil -1}{2},  \text{ if }c=\frac{i}{d}, i=1\ldots d;
$$
\begin{align*}
n_{c,0}(f)=(i-1)(d-i-1)-\sum_{m\ge 3}\nu_m\left(\lceil \frac{im}{d}\rceil -1\right)\left(  m-\lceil \frac{im}{d}\rceil \right),\\
\text{ if }c=\frac{i}{d}+1, i=1\ldots d;
\end{align*}
\begin{align*}
n_{c,0}(f)=\binom{d-i-1}{2} -\sum_{m\ge 3}\nu_m\binom{m-\lceil  \frac{im}{d} \rceil}{2}-\delta_{i,d},\\
\text{ if }c=\frac{i}{d}+2, i=1\ldots d;
\end{align*}
where $
\nu_m=\#\{ P\in D\ |\ mult_PD=m\ \},
$ and $\delta_{i,d}=1$ if $i=d$ and $0$ otherwise.

\me
The motivic, $p$-adic, and topological zeta functions also depend only on the combinatorics, \cite{B-BSY}.

\me
The $b$-function is not a combinatorial invariant, according to a recent announcement of U. Walther. For computations of $b$-functions, by general properties already listed in this survey, it is enough to restrict to the case of so called ``indecomposable central essential" complex arrangements. The {\it $n/d$-Conjecture} says that for such an arrangement of degree $d$, $-n/d$ is a root of the $b$-function. This is known only for reduced arrangements when $n\le 3$, and when $n>3$ for reduced arrangements with $n$ and $d$ coprime and one hyperplane in general position, \cite{B-BSY}. For reduced arrangements as above, it is also known that: if $b_f(-c)=0$ then $c\in (0,2-1/d)$, and $-1$ is a root of multiplicity $n$ of $b_f(s)$, \cite{Sa-barr}.

\me\ni{\it Example.} If $f$ is a generic central hyperplane arrangement, then U. Walther \cite{W}, together with the information about the root $-1$ from above, showed that
$$
b_f(s)=(s+1)^{n-1}\prod_{j=n}^{2d-2} \left (s+\frac{j}{d} \right ).
$$

\me The Monodromy Conjecture holds for all hyperplane arrangements; the Strong Monodromy Conjecture holds for a hyperplane arrangement $D\subset K^n$ if the $n/d$-Conjecture holds, \cite{B-BMT}. In particular, the Strong Monodromy Conjecture holds for all reduced arrangements in $\le 3$ variables, and for the reduced arrangements in $4$ variables of odd degree with one hyperplane in generic position, \cite{B-BSY}.

\subsection {Discriminants of finite reflection groups.} A {\it complex reflection group} is a group $G$, acting on a finite-dimensional complex vector space $V$,  that is generated by elements that fix a hyperplane pointwise, i.e. by complex reflections.  Weyl groups and Coxeter groups are complex reflection groups.   The  ring of invariants is a polynomial ring: $\bC[V]^G=\bC[f_1,\ldots ,f_n]$. Here $n=\dim V $, and
$f_1,\ldots ,f_n$ are some algebraically independent invariant polynomials. The degrees $d_i=\deg f_i$ are determined uniquely. The finite irreducible complex reflection groups are classified by Shephard-Todd, \cite{Br} .

\me
Let $D=\cup_i D_i$ be the union of the reflection hyperplanes, and let $\al_i$ denote a linear form defining $D_i$. Let $e_i$ be the order of the  subgroup fixing $D_i$. Consider the invariant polynomial 
$$
\delta = \prod_i \al_i^{e_i} \in \bC[V]^G.
$$
Viewed as a polynomial in the variables $f_1,\ldots, f_n$, it defines a regular map $\Delta:V/G\cong \bC^n\ra\bC$, called the {\it discriminant}.

\me The $b$-functions of  discriminants of the finite irreducible complex reflection groups have been determined in terms of the degrees $d_i$ for Weyl groups in \cite{O-Weil} and for Coxeter groups in \cite{O-Cox}:
$$
b_\Delta(s)=\prod_{i=1}^n\prod_{j=1}^{d_i-1}\left( s+\frac{1}{2}+\frac{j}{d_i}\right).
$$
In the remaining cases, the zeros of the $b$-functions are determined in \cite{DL-refl}.

\me The monodromy zeta function of $\Delta$ has also been determined in terms of the degrees $d_i$, \cite{DL-refl}.

\subsection {Generic determinantal varieties.}  Let $M$ be the space of all matrices of size $r\times s$, with $r\le s$ . The {\it $k$-th generic determinantal variety} is the subvariety $D^k$ consisting of matrices of rank at most $k$. 

\me The multiplier ideals $\cJ(M,c\cdot D^k)$ have been computed in \cite{Jo}. In particular, the log canonical threshold is
$$
\lct (M,D^k)=\min_{i=0,\ldots ,k}\frac{(r-i)(s-i)}{k+1-i}.
$$

\me The topological zeta function is computed in \cite{Roi}:
$$
Z^{top}_{D^k}(s) =\prod_{c\in\Omega}\frac{1}{1-sc^{-1}}\; ,$$
where
$$
\Omega =\left\{ -\frac{r^2}{k+1}, -\frac{(r-1)^2}{k}, -\frac{(r-2)^2}{k-1},\ldots , -(r-k)^2 \right \}.$$

\me
The number of irreducible components of the $n$-th jet scheme $D^k_n$ is 1 if $k=0, r-1$, and is $n+2-\lceil(n+1)/(k+1)\rceil$ if $0<k<r-1$, \cite{Roi}.

\me It is not known in general how to compute the $b$-function of $D^k$. 

\me\ni{\it Example.} The oldest example of a nontrivial $b$-function is due to Cayley. Let $f=det (x_{ij})$ be the determinant of an $n\times n$ matrix of indeterminates. Then $b_f(s)=(s+1)\ldots (s+n)$ and the differential operator from the definition of the $b$-function is $P=det (\pa /\pa x_{ij})$.
 
\subsection{Prehomogeneous vector spaces.} A {\it  prehomogeneous vector space (pvs)} is a vector space $V$ together with a connected linear algebraic group $G$ with a rational representation $G\ra GL(V)$  such that $V$ has a Zariski dense $G$-orbit. The complement of the dense orbit is called the {\it singular locus}. The pvs is {\it irreducible} if $V$ is an irreducible $G$-module. The pvs is {\it regular} if the singular locus is a hypersurface ${f=0}$. Irreducible regular pvs have been classified into 29 types by Kimura-Sato. The ones in a fixed type are related to each other via a so-called {\it castling transformation}, and within each type there is a "minimal" pvs called {\it reduced}, \cite{Ki}. 

\me
The $b$-functions $b_f(s)$ have been computed for irreducible regular pvs using microlocal calculus by  Kimura (28 types) and Ozeki-Yano (1 type), \cite{Ki}. For an introduction to microlocal calculus see \cite{Ka-book}.

\me The $p$-adic zeta functions of 24 types of irreducible regular pvs have been computed by Igusa. The Strong Monodromy Conjecture has been checked for irreducible regular pvs:  24 types by Igusa \cite{Ig-book, Ki}, and the remaining types by Kimura-Sato-Zhu  \cite{KSZ}.

\me
The castling transform for motivic zeta functions and for Hodge spectrum has been worked out in \cite{Lo-castl} .

\me
There are additional computations of $b$-functions for pvs beyond the case of irreducible and regular ones. We mention a few results. For the reducible pvs, an elementary method
to calculate the $b$-functions of singular loci, which uses the known formula for $b$-functions of one variable, is presented in \cite{Uk}. The decomposition formula for
$b$-functions, which asserts that under certain conditions,
the $b$-functions of reducible pvs have
decompositions correlated to the decomposition of representations, was given in \cite{FS}. By using the decomposition formula, the $b$-functions of relative invariants arising from the quivers of type $A$ have been determined in \cite{Su}. 

\me A {\it linear free divisor} $D\subset V$ is the singular locus of a particular type of pvs. One definition is that the sheaf of vector fields tangent to $D$ is a free $\cO_V$-module and has a basis consisting of vector fields of the type $\sum_jl_{j}\partial_{x_j}$, where $l_j$ are linear forms. Another equivalent definition is that  $D$ is the singular locus of a pvs $(G,V)$ with $\dim G=\dim V=\deg f$, where $f$ is the equation defining $D$. To bridge the two definitions, one has that $G$ is the connected component containing the identity of the group  $\{A\in GL(V)|\ A(D)=D\}$. Quiver representations give often linear free divisors \cite{BM}. The $b$-functions for linear free divisors have been studied in \cite{GSc, Sev} and computed in some cases.

\me\ni{\it Example.} Some interesting examples of $b$-functions, related to quivers of type $A$ and to generic determinantal varieties with blocks of zeros inserted, are computed in \cite{Su}. For example, let $X, Y, Z$ be matrices of three distinct sets of indeterminates of sizes $(n_2,n_1)$, $(n_2,n_3)$, $(n_4,n_3)$, respectively, such that $n_1+n_3=n_2+n_4$ and $n_1< n_2$. Then for
$$
f=\det\left [
\begin{array}{cc}
X & Y\\
0 & Z
\end{array}
\right ],
$$
the $b$-function is 
$$
b_f(s)=(s+1)\ldots (s+n_3)\cdot(s+n_2-n_1+1)\ldots (s+n_2).
$$

\subsection {Quasi-ordinary hypersurface singularities.}  A hypersurface germ $(D,0) \subset (\bC^{n},0)$ is {\it quasi-ordinary} if there exists a finite
morphism $(D, 0) \ra (\bC^{n-1}, 0)$ such that the discriminant locus is contained in a normal crossing divisor. In terms of equations, $f\in\bC[x_1,\ldots ,x_n]$ has quasi-ordinary singularities if the discriminant of $f$ with respect to $y=x_n$ equals $x_1^{u_1}\ldots x_{n-1}^{u_{n-1}}\cdot h$, where $h(0)\ne 0$. 

\me
Quasi-ordinary hypersurface singularities generalize the case of plane curve in the sense that they are higher dimensional singularities with Puiseux expansions. The {\it characteristic exponents} $\lam_1<\ldots <\lam_g$ of an analytically irreducible quasi-ordinary hypersurface germ $f$ are defined as follows. The roots of $f(y)$ are fractional power series $\zeta_i\in\bC[[x_1^{1/e},\ldots ,x_{n-1}^{1/e}]]$, where $e=\deg_y f$. The difference of two roots of $f$ divides the discriminant, hence $\zeta_i -\zeta _j = x^{\lam_{ij}}h_{ij}$, where $h_{ij}$ is a unit. Then $\{\lam_{ij}\}$ is the set of characteristic exponents which we order and relabel it $\{\lam_k\}$. The set of characteristic exponents is an invariant of and equivalent to the embedded topological type of the germ, \cite{Gau}. By a change of variable it can be assumed that the Newton polytope of $f$ is determined canonically by the characteristic exponents. There is a canonical way to relabel the variables and to order the characteristic exponents, \ref{subCanOrd}.

\me There are explicit embedded resolutions of quasi-ordinary singularities in terms of characteristic exponents, \cite{GP}.

\me The monodromy zeta function has been computed in \cite{GP-mon}.

\me The Monodromy Conjecture holds for quasi-ordinary  singularities \cite{AB-mc}.

\me Jet schemes of quasi-ordinary singularities have been analyzed in \cite{Ro, Manuel}. For analytically irreducible quasi-ordinary hypersurface singularities the motivic zeta function, and hence the log canonical threshold, Hodge spectrum and the monodromy zeta function, have been computed in terms of the characteristic exponents in \cite{Manuel}.  

\me\ni{\it Example.}  A refined formula for the log canonical threshold of an analytically irreducible quasi-ordinary hypersurface singularity $f$ was given in \cite{B-qord}. Let $\lam_i$ be ordered as in \ref{subCanOrd}, and let $\lam_{i,j}$ denote the $j$-th coordinate entry of the vector $\lam_i$. Then:

\medskip
(a) $f$ is log canonical if and only if it is smooth, or $g=1$ and
the nonzero coordinates of $\lam_1$ are $1/q$, or $g=1$ and the nonzero coordinates of $\lam_1$ are $1/2$ and  $1$.

\medskip
(b) With $i_1$ and $j_2$ defined as in \ref{subCanOrd}, if $f$ is not log canonical, then: 

\begin{enumerate}
\item if $i_1\ne n-1$, $$\lct (f)=\displaystyle\frac{1+\lam_{1,n-1}}{e\lam_{1,n-1}};$$  

\item if $i_1=n-1, j_2\ne 0$,
 and $\lam_{2,j_2}\ge n_1(\lam_{2,n-1}-\frac{1}{n_1}+1)$, $$\lct (f)=\displaystyle{\frac{1+\lam_{2,j_2}}{\displaystyle\frac{e}{n_1}\cdot\lam_{2,j_2}}}$$ 
 
\item    if $i_1=n-1$ and $j_2=0;$ or  if $i_1=n-1,\; j_2\ne 0$ and $\lam_{2,j_2}< n_1(\lam_{2,n-1}-\frac{1}{n_1}+1)$,
$$\lct (f)=\displaystyle\frac{1+\lam_{2,n-1}}{\displaystyle\frac{e}{n_1}\cdot(\lam_{2,n-1}+1-\displaystyle\frac{1}{n_1})} .$$

\end{enumerate}

\subsubsection{Canonical order.}\label{subCanOrd} 
Let $f$ be an analytically irreducible quasi-ordinary hypersurface. Let
$0=\lam_0<\lam_1<\ldots <\lam_g\in \bZ^{n-1}$ be the
characteristic exponents of $f$. Set for $1\le j\le g$,
\begin{align*}
n_j&:=\# \ (\bZ^{n-1}+\lam_1\bZ+\ldots +\lam_{j}\bZ)/
(\bZ^{n-1}+\lam_1\bZ+\ldots +\lam_{j-1}\bZ).
\end{align*}
It is known that $n_1\ldots n_g=\deg_y f$. 
We set  $n_0=0$. One can permute the variables $x_1,\ldots ,x_n$
such that for $j<j'$ we have $(\lam_{1,j},\ldots ,\lam_{g,j})\le
(\lam_{1,j'},\ldots ,\lam_{g,j'})$ lexicographically.

\medskip

This ordering defines $j_1=n-1>j_2\ge j_3\ge\ldots\ge j_g \ge 0$
such that $j_i=\max\{ j \ \ |\ \lam_{i-1,j}=0  \}$, where we set
$j_i=0$ if $\lam_{i-1,j}\ne 0$ for all $j=1,\ldots n-1$. For $j\le j_i$, $\lam_{i,j}$ can be written as a rational number
with denominator $n_i$. Define $i_k=\max\{  j\le j_k\ |\
\lam_{k,j}=1/n_k \}$ if $j_k>j_{k+1}$, and $i_k=j_{k+1}$ if
$j_k=j_{k+1}$ or $\lam_{k,i}\ne 1/n_k$ for all $i=1,\ldots ,j_k$.
Whenever $j_k<j\le j'\le i_{k-1}$, we have
$\lam_{k,j}\le \lam_{k,j'}$ by our ordering. This defines the notation used in the previous example.

\subsection{Computer programs.}\label{subComp}   Most of the algorithms available for singularity invariants depend on Gr\"{o}bner bases or resolution of singularities. Hence even for an input of small complexity one might not see the end of a computation. Nevertheless, the programs are very useful.

\me As a general rule, one should check the documentation of Singular \cite{Si}, Macaulay2 \cite{M2}, and Risa/Asir \cite{R/A}
for a list of packages pertaining to singularity theory, which can and hopefully does increase often. Singular has the most extensive packages dedicated to singularity theory; see  \cite{GPf}, which contains the theory and subtleties in connection with local computations as opposed to global computations.

\me
Packages involving $\cD$-modules have been implemented extensively. For computation of $b$-functions, T. Oaku gave the first general algorithm. There are now algorithms implemented even for computation of generalized twisted $b$-functions $b_{f,g}(s)$ for ideals $f$ \cite{ABL, Shi, BL}.

\me Using that multiplier ideals can be written in terms of $\cD$-modules \cite{B-Comp}, algorithms are now implemented that compute multiplier ideals and jumping numbers \cite{Shi, BL}.

\me Villamayor's algorithm for resolution of singularities has also been implemented. Some topological zeta functions can be computed using this, \cite{FP}. 

\me For computing the Hodge spectrum and $b$-functions of isolated hypersurface singularities, as well as other invariants, see \cite{Sch}.

\me For nondegenerate polynomials there are programs computing the $p$-adic, topological,  monodromy zeta functions  \cite{Poly}, and the Hodge spectrum \cite{En}.

\me Multiplier ideals for hyperplane arrangements can be computed with \cite{DS-m2}. We are waiting for the authors to implement the combinatorial formulas of \cite{B-SpAr} for jumping numbers and Hodge spectra.

\me Currently there are no feasible algorithms available for computing test ideals or $F$-jumping numbers even for polynomials of very small complexity. This is due to appearance of high powers of the ideals involved in the definitions.

\subsection{Questions.}
\

\me Complete a general additive Thom-Sebastiani formula for $b$-functions. Are there additive Thom-Sebastiani rules for multiplier ideals and test ideals? Generalize \cite{DV} to obtain this rule for the motivic zeta function. This would recover the rule for the motivic Milnor fiber.
 
\me Is there a multiplicative Thom-Sebastiani for $b$-functions?

\me
Let $f$ and $g$ be polynomials. Does there exist a Summation rule that allows one to ``generate" the Verdier specialization functor $\Sp_{(f,g)}$ of the ideal generated by $f$ and $g$ via the nearby cycles $\psi_{f^\lam\cdot g^\mu}$ with $\lam+\mu=1$? Here $\lam$ and $\mu$ would represent eigenvalues of the monodromy along $f$ and $g$, respectively. A similar question can be asked about motivic zeta functions.

\me
Does the Hodge spectrum of a hypersurface germ satisfy a semicontinuity property similar to what happens in the isolated case? 

\me
For reduced and irreducible plane curves, determine the arithmetic invariants needed along with the Puiseux pairs to compute the test ideals and $F$-jumping numbers for a fixed reduction modulo $p$.

\me Prove or correct the formula conjectured by T. Yano for  the $b$-function of a general reduced and irreducible plane curve among those with fixed Puiseux pairs.

\me
Write the Hodge spectrum of nondegenerate polynomials, with non-necessarily isolated singularities, in terms of the Newton polytope.
 
\me Compute the $F$-jumping numbers of hyperplane arrangements, or of some other class of examples besides monomial ideals. 
 
\me
Solve the combinatorial problem that completes the proof of the $n/d$-Conjecture for hyperplane arrangements, and thus of the Strong Monodromy Conjecture for hyperplane arrangements.

\me We have already mentioned the problem of combinatorial invariance of the dimension of the (Hodge pieces of the) cohomology of the Milnor fibers for hyperplane arrangements. This is currently viewed as the ``the holy-grail" in the theory of hyperplane arrangements. The problem fits between the combinatorial invariance of the fundamental group of the complement, which is not true \cite{Ry}, and that of the cohomology of the complement, which is true \cite{OS}.  

\me
What can one say about the other zeta functions, besides the monodromy zeta function, for discriminants of irreducible finite reflection groups? Since the $b$-functions are already determined, maybe the Strong Monodromy Conjecture can be checked.

\me Compute the $b$-function of  generic determinantal varieties. These varieties have certain analogies with monomial ideals, \cite{BV}. Maybe the strategy for computing $b$-functions of monomial ideals can be pushed to work for determinantal varieties.

\me Compute the Hodge spectrum of the 29 types of irreducible regular prehomogeneous vector spaces. By the castling transformation formula of Loeser \cite{Lo-castl}, it is enough to compute the Hodge spectrum for the reduced ones.

\me Construct feasible algorithms for computing test ideals and $F$-jumping numbers. We are also lacking algorithms for Hodge spectra and $p$-adic zeta functions besides the cases mentioned in \ref{subComp}. However, due to their relation with $\cD$-modules and resolution of singularities, it should be possible to give such algorithms.

\end{document}